\documentclass[10pt]{article}
\usepackage{amsmath,amsthm,amsopn,amstext,amscd,amsfonts,amssymb,verbatim}
\usepackage{caption}
\usepackage{subcaption}
\usepackage[dvipsnames]{xcolor}
\colorlet{LightRubineRed}{RubineRed!70}
\colorlet{Mycolor1}{green!10!orange}
\definecolor{Mycolor2}{HTML}{00F9DE}

\usepackage{color}
\usepackage[numbers]{natbib}
\usepackage{graphicx}
\usepackage[explicit]{titlesec}

\pdfoutput=1

\renewenvironment{abstract}
                 {\vspace{6pt}
                  \begin{minipage}{5.3in}
                  \textbf{Abstract--}
                  \noindent\ignorespaces
                 }
                 {\end{minipage} }
                 
\newenvironment{keywords}
                 {\vspace{6pt}
                  \begin{minipage}{5.3in}
                  \textbf{Keywords:}
                  \noindent\ignorespaces
                 }
                 {\end{minipage} }                 
\newtheorem{theorem}{\textbf{Theorem}}[section]

\theoremstyle{definition}
\newtheorem{definition}{\textbf{Definition}}[section]

\titleformat{\section}[block]{\normalfont\center\thesection. #1}{}{0pt}{}
\titleformat{\subsection}[block]{\normalfont\center\thesubsection. #1}{}{0pt}{}

\title{\Large \textbf{Dependent and Independent Time Series}}
\author{
  Fredy O. P\'erez-Ram\'{\i}rez$^{a,b}$, Francisco J. Caro-Lopera$^{c,*}$,\\
  Jos\'e A. D\'{\i}az-Garc\'{\i}a$^{d}$ and Graciela Gonz\'alez-Far\'{\i}as$^{e}$\\[2ex]
  \textit{$^{a}$University of Medellin, Faculty of Engineering,}\\
  \textit{ Medell\'{\i}n, Colombia,} \\
   \textit{$^{b}$University of Medellin, Ph.D. in Modeling and Scientific Computing,}\\
   \textit{Faculty of Basic Sciences, Medell\'{\i}n, Colombia,}\\ 
   \textit{$^{c}$University of Medellin, Faculty of Basic Sciences,}\\ 
   \textit{ Medell\'{\i}n, Colombia,} \\
   \textit{$^{b}$Universidad Aut\'onoma de Chihuahua, Facultad de Zootecnia y Ecolog\'{\i}a,}\\
   \textit{ Chihuahua, M\'exico}\\
    \textit{Centro de Investigación en Matemáticas,}\\
    \textit{Departamento de Probabilidad y Estadística,Guanajuato, M\'exico}\\
    \textit{foperez@udemedellin.edu.co, fjcaro@udemedellin.edu.co},\\ \textit{jadiaz@cimat.mx, farias@cimat.mx}.
    ${}^{*}$ Corresponding Author.
}
\begin{document}
\maketitle

\begin{abstract}
The time series theory is set in this work under the domain of general elliptically contoured distributions. The advent of a time series approach that is in accordance with the expected reality of dependence between errors, transfers the increasingly complex and difficult to handle correlation analysis into a discipline that models volatility from a new view of a likelihood based on dependent probabilistic samples. The equally important problem of model selection is strengthened, but at the same time criticized with the introduction of degrees of evidence of significant difference in the modified BIC criterion . The demanding scale of differentiation puts a well-known database in trouble by observing insignificant relevance between the hierarchical models most used in the theory of time series under independence, such as Arch, Garch, Tgarch and Egarch. For extreme cases where the probabilistic independence of the samples is exceptionally demonstrated by an expert, the article also proposes the theory of time series under elliptical models, but with the same demanding comparison of the degrees of evidence of differences of the modified BIC. The example studied under this approach also does not denote any advantage of the hierarchical models studied. Such a new perspective for a likelihood based on dependent probabilistic samples has arisen naturally in similar context in finance under the setting of multivector variate distributions.    

\end{abstract}



\section{Introduction}\label{sec:introduction}
Time series theory has been widely explored over the past four decades. Among other issues, volatility modelling has led to the formulation of several theories that improve upon previous theories by eliminating major assumptions.
We focus on univariate models for volatile time series. It is of interest to model phenomena where the conditional variance is changing over time, a very common situation in the financial sector. There, the investor is interested in making predictions about the rate of return and its volatility while holding a certain stock. Anticipating performances and risks over short and long time periods are critical for both the investor and the issuer. Both are interested in the mean return, but the former focuses on the short term, revealed in the conditional variance, while the latter is interested in the unconditional variances for long periods. Time series associated with stock prices, exchange rates, monetary aggregates and inflation rates generally do not have constant averages and are governed by high volatilities, that is, large fluctuations in variability with respect to the mean. 
This characteristic, deeply rooted in finance and economics, has given rise to a large number of theories and research explaining volatility. The initial works are due to \citet{Engle1982}, with the ARCH models (Autoregressive Conditional Heterocedasticity model) and to \citet{Bollerslev1986}, with the GARCH models (generalised ARCH). GARCH is understood as the Generalized Autoregressive Conditional Heterocedasticity model; heteroskedasticity is taken as the variation in variance over time. The term conditional implies a dependence on observations from the immediate past and autoregressive describes the feedback mechanism, which incorporates past observations in the present. GARCH then, is a mechanism that includes past variances in the explanation of future variances. More specifically, GARCH is a time series model that uses past variances and predictions of past variances to predict future variances. \citet{EngleBollerslev1986} with the GARCH-M model. Other approaches include Threshold GARCH (TGARCH) of \citet{Zakoian1994}; GJR GARCH of \citet{Glostenetal1993}; EGARCH by  \citet{Nelson1991}; \citet{DingGrangerEngle1993}  with the APARCH model; \citet{BaillieBollerslevMikkelsen1996} with the FIGARCH model; among many others.

The models have evolved from the weaknesses of their predecessors. The theoretical and practical evolution of each of the models is certainly intricate, revealing the criticisms of the predecessors and attempting to improve the explanation of the elusive volatility. Regardless of the complexity of each, there are several immutable underlying elements that make them similar. This immutability seems not to be disputed or considered, as long as there is a certain historical consensus crossed by the ease of interpretation and the reduction of computational cost. Among others, two fundamental elements are not sufficiently explained in the addressed time series models: the grades of evidence in the model selection criteria and the use of probabilistically independent samples. 

First, for our proposal, we will focus on the Yang and Yang criterion (\citet{YangYang2007}), which in turn is aligned with a strict selection criterion based on the degree of evidence of difference. Degrees of evidence that, when applied to time series, represent a paradigm shift in some examples such as those presented here. The modification of the BIC statistics arrived in \citet{YangYang2007} after its conceptualization in a coding theory environment given by \citet{Rissanen1978}.   
The modified BIC$^{*}$ is defined as follows:
\begin{equation}\label{BIC*}
    \mbox{BIC$^{*}$}=\left[-2\mathcal{M}+n_{p}(\log(n+2)-\log 24)\right]/n,
\end{equation}
where $\mathcal{M}$ is the maximum log-likelihood, $n$ is the sample size and $n_{p}$ is the number of estimated parameters. 
Now, for a comparison of two models, \citet{KassRaftery1995} and \citet{Raftery1995} proposed the grades of evidence of Table \ref{Table:BIC*} for BIC$^{*}$ difference.
\begin{table}[h!]
\centering
\begin{tabular}{||c c||} 
 \hline
 BIC$^{*}$ difference & Evidence \\ [0.5ex] 
 \hline\hline
 0-2 & Weak \\ 
 2-6 & Positive \\
 6-10 & Strong \\
 $>$10& Very strong \\ [1ex] 
 \hline
\end{tabular}
\caption{Grades of evidence for BIC$^{*}$ difference}
\label{Table:BIC*}
\end{table}
Then, we have reached a crucial point in the proposal of this article, since the degrees of evidence for a very strong difference between two time series models assume values that are difficult to reach from the theoretical improvements that the theories are implementing.
This particular work exhibits a counterexample in which none of the models reach even a positive difference with the simplest model. Between Arch(1) and any Garch(p,q), the degree of evidence does not reach 2. This aspect is contrary to what the theory indicates.
Similar conclusions are obtained between the nested models Tgarch and Egarch, when comparing them with their respective simplest model.  Some expensive commercial time series software, books and articles rank models from minimum differences of thousandths or less in the traditional Akaike (AIC), Schwarz (BIC) or Hannan-Quinn (HQC) criteria. An analysis with the Yang and Yang criterion from the indicated degrees of evidence leaves these classical decisions with the minimum evidence of weak. The lack of incorporation of degrees of evidence between time series models and their importance is detailed due to a major fact that does not allow to accurately show the effectiveness of a particular model. It is about the difficulty to compare volatility models from the estimation of the variance that they determine. Unlike the mean, which can be compared between models, and is completely evident and comparable with the real and adjusted data, the observed volatility cannot be compared directly with the estimates, so there is no criterion of reality and field that tells us which is the best criterion, then functioning as a simple resource, the RMSE or the minimum difference, no matter how much between a higher and a lower criterion. The RMSE, used for comparing non-nested models, represents the same uncertainty noted to differentiate models by thousandths in the change of the AIC; BIC or HQC, because their distribution is not known and the respective probabilistic evidence cannot be calculated. So again evidence grades of \citet{KassRaftery1995} and \citet{Raftery1995} are of great value to this research.

Thus, the inclusion of degrees of evidence of difference in the Yang-Yang criterion represents an important paradigm shift that must be considered in time series models under independent probability samples.

We come to another crucial point discussed in this paper, which is the classical estimation of the parameters of simple or complex volatility models, from classical likelihood. A function underlying a principle that is difficult to assume on independent probability samples.

The concept of likelihood has remained unchanged in time series theory, leaving the deeper problem of the real dependence of samples on time to the correction or study of correlation, which in even more extreme cases results in tests that fail the expected correlation between events that are actually related in time.

Imposing probabilistic independence on the complex time series of postmodern markets, highly influenced by the changing social, political, cultural and economic conditions of the world, results in a contradiction to the common sense that the investor knows in the empirical knowledge and very broad conception of the evident sequential interaction of the economic variables.

Including the explicit reality of the change in financial series over time means considering models of probabilistic dependence in the samples. To justify the incorporation of probabilistic dependence in the samples, as a fundamental element in the modelling of time series, it is enough to observe the example in \citet{dgclpr:22} and its impact on the decision in finance. There we can see how a traditional decision based on probabilistic independence appears in the tails of the real distribution based on dependence, an aspect that denotes as a strange event, that which is supposed to be abundant in the mean of the process under consideration. 

An alternative to deal with the problem of probabilistic dependence in time series can be considered copula theory. However, the estimation of the parameters maintains the use of classical likelihood. Artificial intelligence, for example neural networks, has the same difficulty when assuming activation functions that only seek to maintain the domain of the events they train. Since they do not work with catastrophic events or small, never-trained samples, their estimates are still off the expected forecast. Even at the height of artificial intelligence's advent, powerful artificial intelligence techniques failed to predict the recent fall of August 5, 2024 on the Tokyo stock exchange. And although warned by the Tokyo catastrophe due to the time difference, AI algorithms also failed to predict what would happen on the New York stock exchange hours later. The paradigm of over- or underestimation in neural networks places activation functions in the same role as the normal distribution for regression models. Back propagation in artificial intelligence processes also end up in the same global minimum search problems, relegating their effectiveness only to local minima. 

Given such difficulties, it is still plausible to consider new scenarios from a mathematical perspective, which at least reduce the latent problems of classical time series models. In particular, redefining a likelihood based on dependent probabilistic samples, which also allows for the involvement of distribution classes, can become a theoretical alternative, which is currently lacking in artificial intelligence applied to time series.

This discussion is placed in Section \ref{sec:main} by setting the time series into the class of elliptically contoured models for describing the volatility in the more realistic probabilistic dependent samples. The independent case under the referred family of distribution is also derived. The volatility models under consideration include Arch, Garch, Tgarch and Egarch. Section \ref{Applications} studies the well known data base of the Gillette time series returns under independent probabilistic samples based on the Kotz
distribution and dependent probabilistic samples based on the Pearson type VII
distribution 

\section{Main result}\label{sec:main}
A new setting for time series based on probabilistic dependent samples require a supporting class of distributions for the random variable. We then introduce the theory of elliptically contoured vector variate distributions, derived from the matrix variate case:
\begin{definition}\label{def1}
Let $\mathbf{X}$ be a $m\times n$ dimensional random matrix whose distribution is
absolutely continuous. Then, the elliptically contoured p.d.f. of the matrix $\mathbf{X}$, denoted by $\mathbf{X}\sim E_{m,n}(\mathbf{M},\boldsymbol{\Sigma} \otimes\boldsymbol{\Phi},h)$, 
has the form
\begin{equation*}
f (\mathbf{X}) = |\boldsymbol{\Sigma}|^{-n/2}|\boldsymbol{\Phi}|^{-m/2}g(tr((\mathbf{X}-\mathbf{M})'\boldsymbol{\Sigma}^{-1}(\mathbf{X}-\mathbf{M})\boldsymbol{\Phi}^{-1})),    
\end{equation*}
where $\boldsymbol{\Sigma}$ and $\boldsymbol{\Phi}$ are positive definite $m\times m$ and $n\times n$ matrices, respectively, and the generator function $g$ satisfies that $0 <\int_{0}^{\infty}v^{mn-1}g(v^{2})dr<\infty$.

The probabilistic dependent samples will required the random vector case. It is obtained when $n=1$. Next, it  is listed throughout specific models depending of simplified variable $v=(\mathbf{x}-\mathbf{\mu})'\boldsymbol{\Sigma}^{-1}(\mathbf{x}-\mathbf{\mu})$ which can be applied as supporting distributions of the time series error law. Here $\mathbf{x}, \boldsymbol{\mu}$ are $m-$dimensional vectors and $\boldsymbol{\Sigma}$ is a positive definite $m\times m$ matrix.   
\begin{itemize}
    \item Kotz Type Distribution: $f(v)=\frac{sr^{(2Q+m-2)/(2s)}\Gamma(m/2)|\boldsymbol{\Sigma}|^{-1/2}}{\pi^{m/2}\Gamma\left((2Q+m-2)/(2s)\right)}v^{Q-1}\exp{-rv^{s}}$, where $r>0, s>0, 2Q+m>2$. If $r=1/2, s=Q=1$, the multivariate normal is obtained. 
    \item Pearson Type VII Distribution: $f(v)=\frac{\Gamma(Q)|\boldsymbol{\Sigma}|^{-1/2}}{(r\pi)^{m/2}\Gamma\left((2Q-m)/2\right)}\left(1+\frac{v}{r}\right)^{-Q}$, where $r>0, Q>m/2$. The multivariate t of $r$ degrees of freedom is obtained when $Q=(m+r)/2$. 
    \item Pearson Type II Distribution: $f(v)=\frac{\Gamma\left(Q+1+m/2\right)|\boldsymbol{\Sigma}|^{-1/2}}{\pi^{m/2}\Gamma\left(Q+1\right)}\left(1-v\right)^{Q}$, where $Q>-1, v\leq 1$.  
    \item Bessel Distribution
    $f(v)=\frac{2^{-Q-m+1}|\boldsymbol{\Sigma}|^{-1/2}v^{Q/2}}{\pi^{m/2}r^{m+Q}\Gamma\left(Q+m/2\right)}K_{Q}\left(\frac{v^{1/2}}{r}\right)$, where $Q$ is integer, $K_{Q}(z)=\pi(I_{-Q}(z)+I_{Q}(z))/(2\sin(Q\pi))$, $|\arg(z)|<\pi$, $|z|<\infty$ and $I_{Q}(z)=\sum_{k=0}^{\infty}\frac{(z/2)^{Q+2k}}{k!\Gamma(k+Q+1)}$. The multivariate Laplace distribution is obtained when $Q=0, r=s/\sqrt{2}, s>0$.
\end{itemize}
Finally, for applications in probabilistic independent samples, we will required the univariate case, which is reached when $m=1$. 
\end{definition}

The models for volatility that we will research in this paper are the most frequently used, namely, Arch, Garch, Tgarch and Egarch. 
As usual, let $y_{t}$ be a random variable measured at time $t$ with an available information on the  $t-1$ past records. Based on the past time, $y_{t}$ in the univariate and the new vectorial way, that we will consider below, are ruled by an elliptical law with mean and variance. Specifically, the error has the form  $\varepsilon_{t}=y_{t}-(1,x_{1,t},\ldots,x_{k,t})(\beta_{0},\beta_{1},\ldots,\beta_{k})'$, where $\beta_{0},\beta_{1},\ldots,\beta_{k}$ are the mean parameters to be estimated jointly with corresponding parameters of the volatility model implicit in the variance $\sigma_{t}^{2}$.

Now, assume that $\mathbf{y}\sim E_{m}(\boldsymbol{\mu},\boldsymbol{\Sigma},h)$ has a multivariate elliptically contoured distribution with probability density function, respect the Lebesgue measure, of parameters $\boldsymbol{\mu}\in \Re^{m}$ and $\boldsymbol{\Sigma}>0$, thus
\begin{equation}\label{multivariateelliptical}
    f_{\mathbf{y}}(\mathbf{y})=c_{m}|\boldsymbol{\Sigma}|^{-1/2}h\left[(\mathbf{y}-\mathbf{\mu})'\boldsymbol{\Sigma}^{-1}(\mathbf{y}-\mathbf{\mu})\right],
\end{equation}
with $\mathbf{y}\in \Re^{m}$, $\int_{0}^{\infty}v^{m-1}h(v^{2})dv<\infty$ and $c_{m}$ is the normalization constant, such that $\int_{\Re ^{m}} f_{\mathbf{y}}(\mathbf{y})d\mathbf{y}=1$.

The univariate case, when $m=1$, takes the form:
\begin{equation}\label{univariateelliptical}
    f_{y}(y)=c_{1}|\boldsymbol{\Sigma}|^{-1/2}h\left[\frac{(y-\mu)^2}{\sigma^2}\right],
\end{equation}
with $y\in \Re$ and $\int_{0}^{\infty}h(v^{2})dv<\infty$.

Multivariate and univariate particular cases can be seen in Definition \ref{def1}, then they induced the following main result:

\begin{theorem}\label{MainTh}
    Set $\boldsymbol{\varepsilon}=\left(\varepsilon_{1},\varepsilon_{2},\ldots,\varepsilon_{T}\right)' \sim E_{T}\left(\mathbf{0},diag\left(\sigma_{1}^{2},\sigma_{2}^{2},\ldots,\sigma_{T}^{2},\right),h\right)$. 
Thus the likelihood of $\boldsymbol{\varepsilon}$ representing a probabilistic dependent sample is given by:
\begin{equation}\label{likelihooddependent}    L\left(\boldsymbol{\varepsilon}|\sigma_{1}^{2},\ldots,\sigma_{T}^{2}\right)=\frac{c_{T}}{\left(\displaystyle\prod_{t=1}^{T}\sigma_{t}^{2}\right)^{1/2}}h\left(\displaystyle\sum_{t=1}^{T}\frac{\varepsilon_{t}^{2}}{\sigma_{t}^{2}}\right).
\end{equation}
Meanwhile the classical likelihood based on a probabilistic independent sample is obtained as:
\begin{equation}\label{likelihoodindependent}    L\left(\varepsilon_{1},\varepsilon_{2},\ldots,\varepsilon_{T}|\sigma_{1}^{2},\ldots,\sigma_{T}^{2}\right)=\frac{c_{1}^{T}}{\left(\displaystyle\prod_{t=1}^{T}\sigma_{t}^{2}\right)^{1/2}}\displaystyle\prod_{t=1}^{T}h\left(\frac{\varepsilon_{t}^{2}}{\sigma_{t}^{2}}\right).
\end{equation}
And given an Ar(k) law, with parameters $\beta_{i}, i=1,\ldots,k$, for the mean, the volatility model for the variance $\sigma_{t}^{2}, t=1, \ldots,T$ takes the following forms:
\begin{itemize}
    \item Autoregressive Conditional Heterocedasticity, Arch(p) (\citet{Engle1982}): 
The loglikelihood of (\ref{likelihooddependent}) and (\ref{likelihoodindependent}) are maximized in terms of the $k+1+p+1$ parameters $\beta_{0},\ldots,\beta_{k},\alpha_{0},\ldots,\alpha_{p}$, where
$$\sigma_{t}^{2}=\alpha_{0}+\sum_{i=1}^{p}\alpha_{i}\varepsilon_{t-i}^{2}.$$ 
    \item Generalized Autoregressive Conditional Heterocedasticity, Garch(p,q) (\citet{Bollerslev1986}): 
In this case the loglikelihood of (\ref{likelihooddependent}) and (\ref{likelihoodindependent}) are maximized in terms of the $k+1+p+1+q$ parameters $\beta_{0},\ldots,\beta_{k},\alpha_{0},\ldots,\alpha_{p},\gamma_{1},\ldots,\gamma_{q}$,
where 
$$\sigma_{t}^{2}=\alpha_{0}+\sum_{i=1}^{p}\alpha_{i}\varepsilon_{t-i}^{2}+\sum_{i=1}^{q}\gamma_{i}\sigma_{t-i}^{2}.$$
    \item Threshold Garch, Tgarch(p,q) (\citet{Zakoian1994}):
In this case the loglikelihood of (\ref{likelihooddependent}) and (\ref{likelihoodindependent}) are maximized in terms of the $k+1+p+1+q$ parameters 
$\beta_{0},\ldots,\beta_{k},\alpha_{0},\ldots,\alpha_{p},$

$\gamma_{1},\ldots,\gamma_{q},\delta_{1},\ldots,\delta_{p}$, where
$$\sigma_{t}^{2}=\alpha_{0}+\sum_{i=1}^{p}\alpha_{i}\varepsilon_{t-i}^{2}+\sum_{i=1}^{q}\gamma_{i}\sigma_{t-i}^{2}+\sum_{i=1}^{p}\delta_{i}N_{t-i}\varepsilon_{t-i}^{2},$$
with $N_{t-i}$ is 1 or 0 if $\varepsilon_{t-i}<0$ or $\varepsilon_{t-i}\geq 0$, respectively.
    \item Exponential Generalized Autoregressive Conditional Heteroscedastic, Egarch(p,q) (\citet{Nelson1991}): 
    The loglikelihood of (\ref{likelihooddependent}) and (\ref{likelihoodindependent}) are maximized in terms of the $k+1+p+1+q+p$ parameters
    $\beta_{0},\ldots,\beta_{k},\alpha_{0},\ldots,\alpha_{p},\gamma_{1},\ldots,\gamma_{q},\delta_{1},\ldots,\delta_{p}$,
where the volatility is defined in the same size parameters of Tgarch as
$$\sigma_{t}^{2}=\exp\left[\alpha_{0}+\sum_{i=1}^{p}\alpha_{i}\left(|\varepsilon_{t-i}/\sigma_{t-i}|-E(Q)\right)+\sum_{i=1}^{q}\gamma_{i}\varepsilon_{t-i}/\sigma_{t-i}+\sum_{i=1}^{p}\delta_{i}\log\sigma_{t-i}^{2} \right]$$
and $E(h)$ depends on the elliptical model $h$. 
\end{itemize}
\end{theorem}

Theorem \ref{MainTh} can be applied now for a number of situations. In particular, modeling volatilities is of great interest in time series, then we focuos on such approach. The independent case is exemplified with the Kotz distribution and the dependent setting is considered under a Pearson Type VII model.

\subsection{Independent probabilistic samples based on the Kotz distribution}\label{subsec:IndependentKotz}

Taking the generator function $h(y)=y^{Q-1}\exp(-y/2)$, as a conditional density based on past information, we have the  loglikelihood  (\ref{likelihoodindependent}) for jointly estimation of Kotz shape parameter $Q>1/2$, and the corresponding parameters in the mean and variance of Arch(p), Garch(p,q), Tgarch(p,q) and Egarch(p,q). In all cases $y_{t}$ is the random variable measured at time $t$ and $\varepsilon_{t}=y_{t}-(1,x_{1,t},\ldots,x_{k,t})(\beta_{0},\beta_{1},\ldots,\beta_{k})'$, where $\beta_{0},\beta_{1},\ldots,\beta_{k}$ are the mean parameters to be estimated jointly with corresponding parameters of the variance $\sigma_{t}^{2}$. The claimed loglikelihood takes de form:
$$
logL\left(\mbox parameters|\varepsilon_{1},\ldots,\varepsilon_{T},\sigma_{1}^{2},\ldots,\sigma_{T}^{2}\right)=T\left[(1/2-Q)\log(2)-\log\Gamma(-1/2+Q)\right]\hspace{5cm}
$$
\begin{equation}\label{loglikelihoodKotzIndependent}    
+\sum_{t=1}^{T}\left[(-1/2)\log(\sigma_{t}^{2})-(1/2)\varepsilon_{t}^{2}/\sigma_{t}^{2}+(Q-1)(\log(\varepsilon_{t}^{2})-log(\sigma_{t}^{2})) \right],
\end{equation}
Next section will provide an application of the four models with a real data set of time series. 

\subsection{Dependent probabilistic samples based on the Pearson type VII distribution}\label{subsec:DependentPearsonVII}
We explore the main model of the paper about the dependent time series, by studying the generator function $f(y)=\left(1+\frac{y}{r}\right)^{-Q}$,where $r>0, Q>m/2$ and $Q=(m+r)/2$. Thus the multivariate t of $r$ degrees of freedom emerges as the conditional density for past information. In this setting, for a given degree of freedom $r$, the loglikelihood based on (\ref{likelihooddependent}) is considered for estimation of the parameters in the mean and variance under the models Arch(p), Garch(p,q), Tgarch(p,q) and Egarch(p,q). In all cases $y_{t}$ is the random variable measured at time $t$ and $\varepsilon_{t}=y_{t}-(1,x_{1,t},\ldots,x_{k,t})(\beta_{0},\beta_{1},\ldots,\beta_{k})'$, where $\beta_{0},\beta_{1},\ldots,\beta_{k}$ are the mean parameters to be estimated jointly with corresponding parameters of the variance $\sigma_{t}^{2}$. The claimed loglikelihood takes de form:
$$
logL\left(\mbox parameters|\varepsilon_{1},\ldots,\varepsilon_{T},\sigma_{1}^{2},\ldots,\sigma_{T}^{2}\right)=\log\left[\Gamma\left((T+r)/2\right)\right]-(T/2)\log\left[\pi r\right]-\log[\Gamma(r/2)]\hspace{5cm}
$$
\begin{equation}\label{loglikelihoodPearsonDependent}    
-(1/2)\sum_{t=1}^{T}\log(\sigma_{t}^{2})-((T+r)/2)\log\left[1+(1/r)\sum_{t=1}^{T}\varepsilon_{t}^{2}/\sigma_{t}^{2}\right],
\end{equation}
The models for volatility are the same of the independent case without considering the Kotz parameter $Q$. Under the Pearson Type VII model we cannot estimate the degrees of freedom $r$, only we consider the corresponding parameters $\beta_{0},\ldots,\beta_{k},\alpha_{0},\ldots,\alpha_{p},\gamma_{1},\ldots,\gamma_{q},\delta_{1},\ldots,\delta_{p}$ in the models Arch(p), Garch(p,q), Tgarch(p,q) and Egarch(p,q).
Finally, we will also explore the preceding models under a real time series data set.

\section{Applications}\label{Applications}
We now consider the application of Theorem \ref{MainTh} into the setting of subsections \ref{subsec:IndependentKotz} and \ref{subsec:DependentPearsonVII}. 

\subsection{Gillette time series returns under independent probabilistic samples based on the Kotz distribution}\label{subsec:GuilletteIndependentKotz}

\begin{figure}
     \centering
     \begin{subfigure}[b]{0.45\textwidth}
         \centering
         \includegraphics[width=\textwidth]{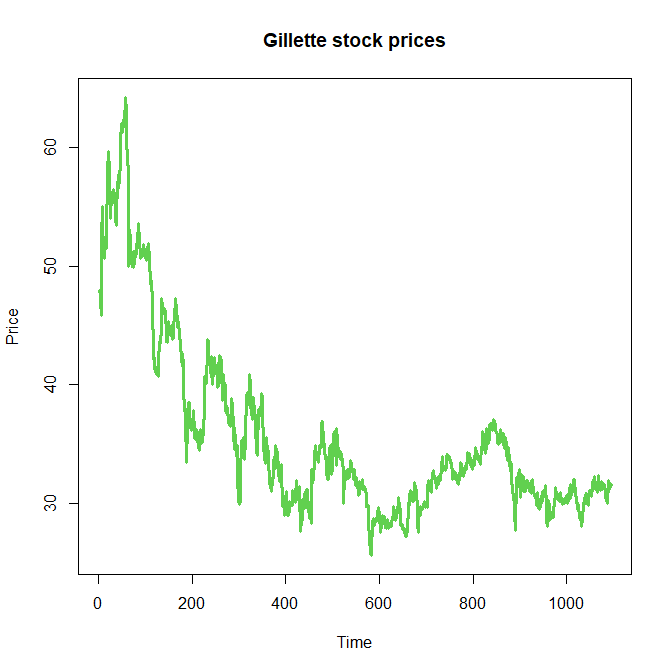}
         \caption{Gillette stock prices}
         \label{prices}
     \end{subfigure}
     \hfill
     \begin{subfigure}[b]{0.45\textwidth}
         \centering
         \includegraphics[width=\textwidth]{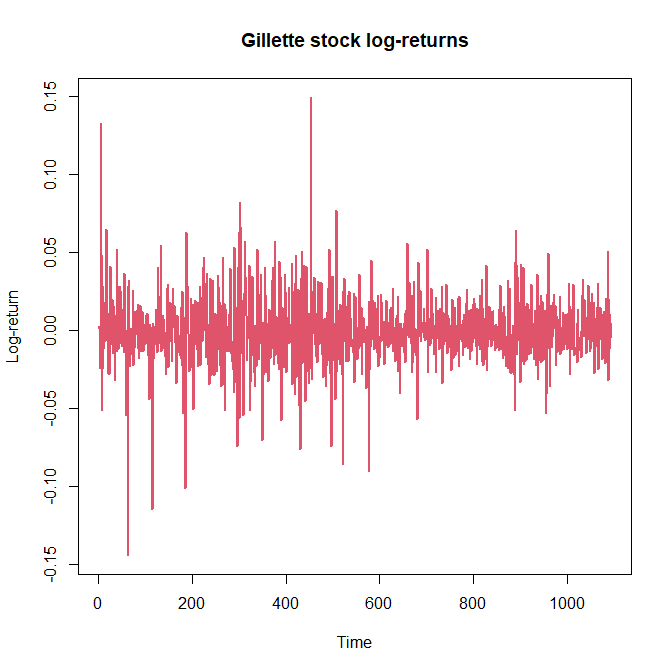}
         \caption{Gillette stock returns}
         \label{logreturns}
     \end{subfigure}
        \caption{Gillette stock time series from January 4-1999 to May 13-2003}
       \label{figurespriceslogreturns}
\end{figure}

Consider the $T=1095$ Gillette stock prices ($p_{t}$, t=1,\ldots T) from January 4-1999 to May 13-2003 (Figure \ref{prices}). In terms of the continuously
compounded return or log return $r_{t}=p_{t}-p_{t-1}$, Figure \ref{logreturns} exhibit a high volatility or conditional variance of the series. Despite the available models for volatility (Arch, Garch, Tgarchm Egarch, etc..), a crucial hidden weak of all time series estimation resides in the use of a likelihood based on an unrealistic independent probabilistic sample. The underlying errors of the assest return models has been considered under different laws, as Gaussian, generalized error, among others. In this subsection we consider the Kotz model of parameter $Q>1/2$, containing the normal case when $Q=1$. 

A optimization routine with the Optimx Package Software in R, provides the maximum loglikelihood estimation of the parameters $Q,\beta_{0},\ldots,\beta_{k},\alpha_{0},\ldots,\alpha_{p},\gamma_{1},\ldots,\gamma_{q},\delta_{1},\ldots,\delta_{p}$ in the models Arch(p), Garch(p,q), Tgarch(p,q) and Egarch(p,q).

\begin{figure}
     \centering
     \begin{subfigure}[b]{0.45\textwidth}
         \centering
         \includegraphics[width=\textwidth]{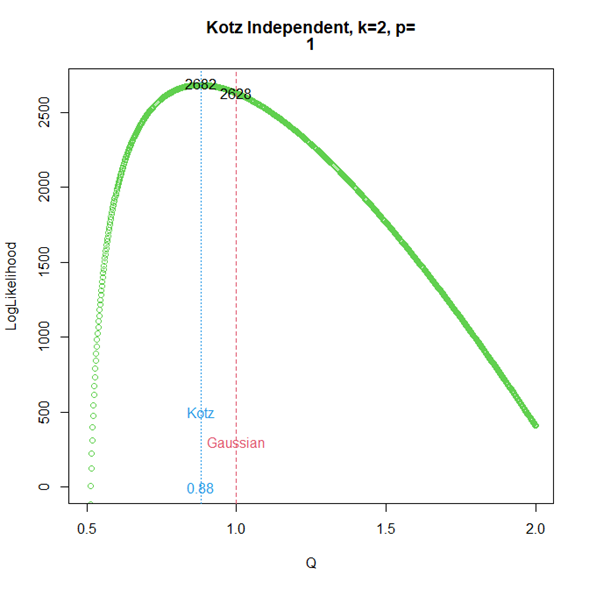}
         \caption{Loglikelihood profile for $Q$.}
         \label{ArchKotzk2p1f1}
     \end{subfigure}
     \hfill
     \begin{subfigure}[b]{0.45\textwidth}
         \centering
         \includegraphics[width=\textwidth]{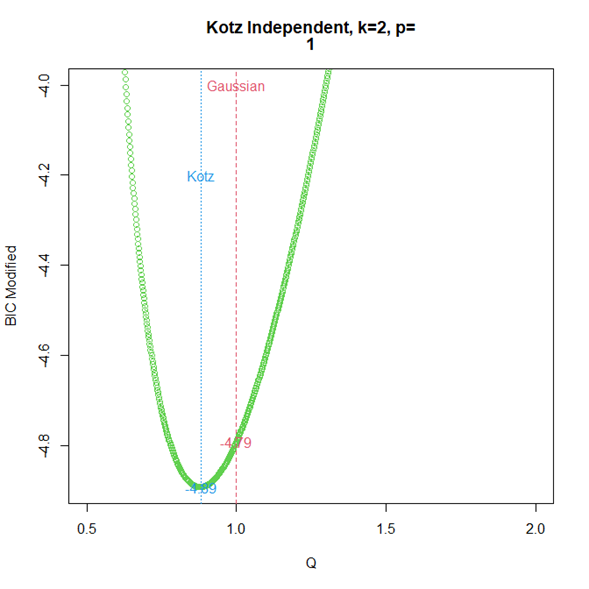}
         \caption{Yang and Yang profile criterion for $Q$.}
         \label{ArchKotzk2p1f2}
     \end{subfigure}
        \caption{Arch Kotz Independent ML estimation for $p=1$.}
       \label{ArchKotzIndependent}
\end{figure}

The profile loglikelihood for the parameter $Q$ in the Arch($p=1$) Kotz model is given in Figure \ref{ArchKotzk2p1f1}, and the parallel behavior with the modified BIC criterion of Yang and Yang, can be seen in Figure \ref{ArchKotzk2p1f2}.

\begin{figure}
     \centering
     \begin{subfigure}[b]{0.45\textwidth}
         \centering
         \includegraphics[width=\textwidth]{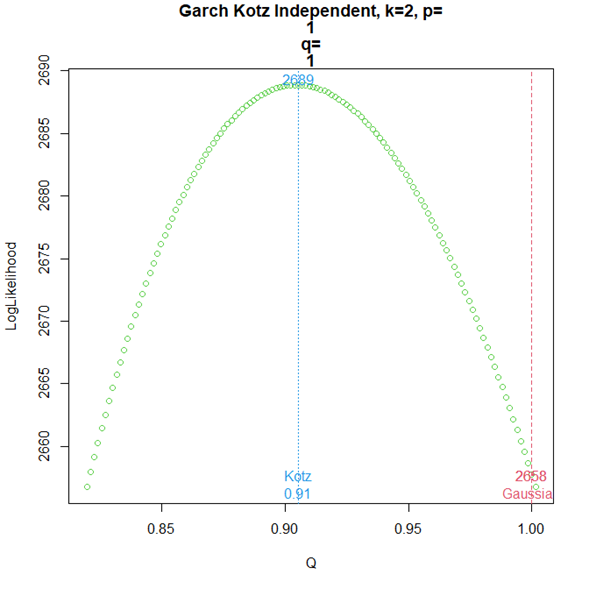}
         \caption{Loglikelihood profile for $Q$.}
         \label{GarchKotzk2p1q1f1}
     \end{subfigure}
     \hfill
     \begin{subfigure}[b]{0.45\textwidth}
         \centering
         \includegraphics[width=\textwidth]{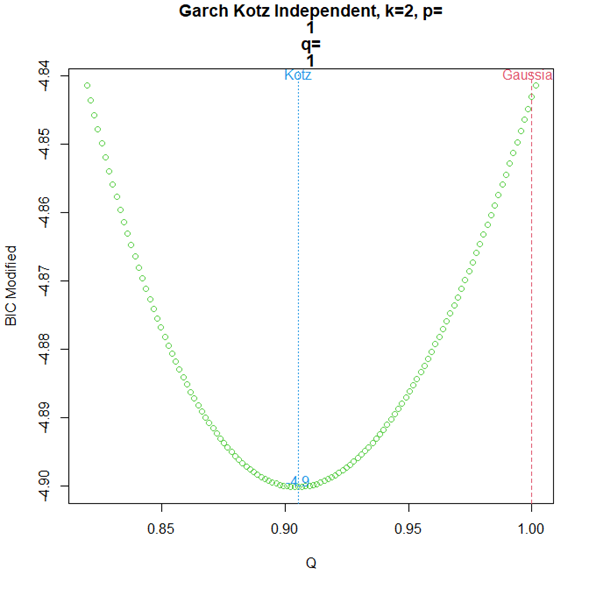}
         \caption{Yang and Yang profile criterion for $Q$.}
         \label{GarchKotzk2p1q1f2}
     \end{subfigure}
        \caption{Garch Kotz Independent ML estimation for $p=1, q=1$.}
       \label{GarchKotzIndependent}
\end{figure}

Figure \ref{GarchKotzk2p1q1f1} shows the profile loglikelihood for the parameter $Q$ in the Garch($p=1$,$q=1$) Kotz independent model. The corresponding modified BIC criterion of Yang and Yang for the profile behavior of $Q$ is given in Figure \ref{GarchKotzk2p1q1f2}.

\begin{figure}
     \centering
     \begin{subfigure}[b]{0.45\textwidth}
         \centering
         \includegraphics[width=\textwidth]{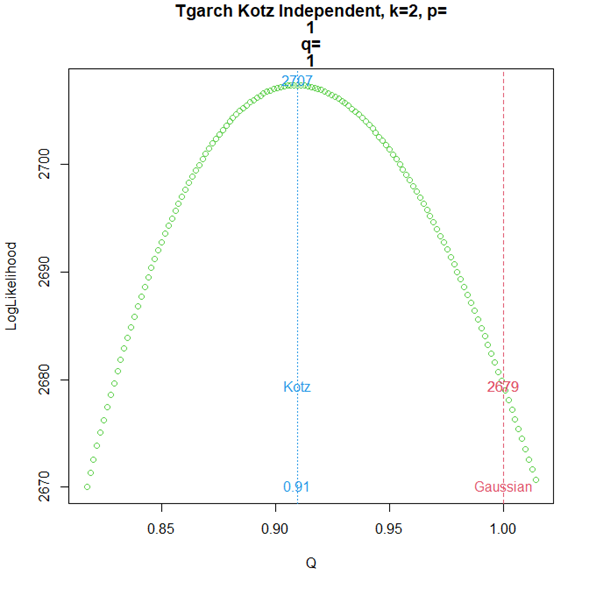}
         \caption{Loglikelihood profile for $Q$.}
         \label{TgarchKotzk2p1q1f1}
     \end{subfigure}
     \hfill
     \begin{subfigure}[b]{0.45\textwidth}
         \centering
         \includegraphics[width=\textwidth]{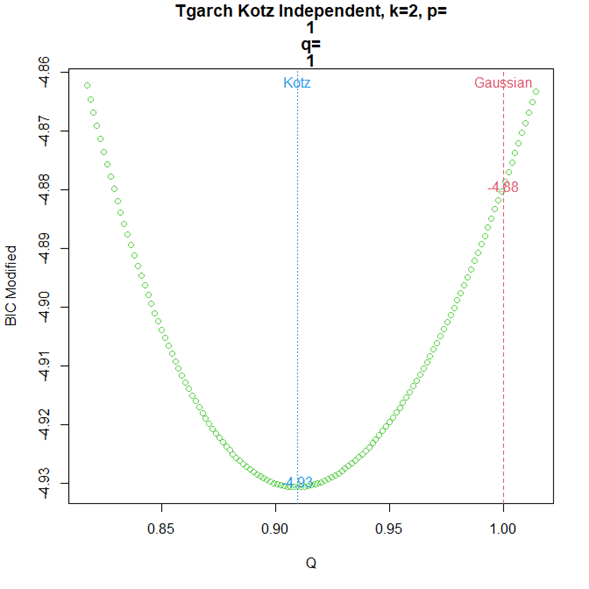}
         \caption{Yang and Yang profile criterion for $Q$.}
         \label{TgarchKotzk2p1q1f2}
     \end{subfigure}
        \caption{Tgarch Kotz Independent ML estimation for $p=1, q=1$.}
       \label{TgarchKotzIndependent}
\end{figure}

In a similar way a Tgarch($p=1, q=1$) Kotz model apparently favors a non Gaussian model as is shown Figure \ref{TgarchKotzk2p1q1f1} for the profile loglikelihood of parameter $Q$. Meanwhile Figure \ref{TgarchKotzk2p1q1f2} exhibit the same conclusion from the modified BIC criterion of Yang and Yang.

\begin{figure}
     \centering
     \begin{subfigure}[b]{0.45\textwidth}
         \centering
         \includegraphics[width=\textwidth]{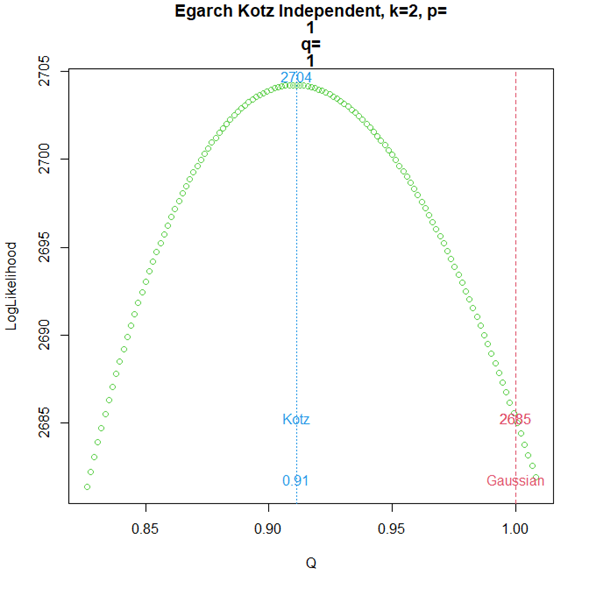}
         \caption{Loglikelihood profile for $Q$.}
         \label{EgarchKotzk2p1q1f1}
     \end{subfigure}
     \hfill
     \begin{subfigure}[b]{0.45\textwidth}
         \centering
         \includegraphics[width=\textwidth]{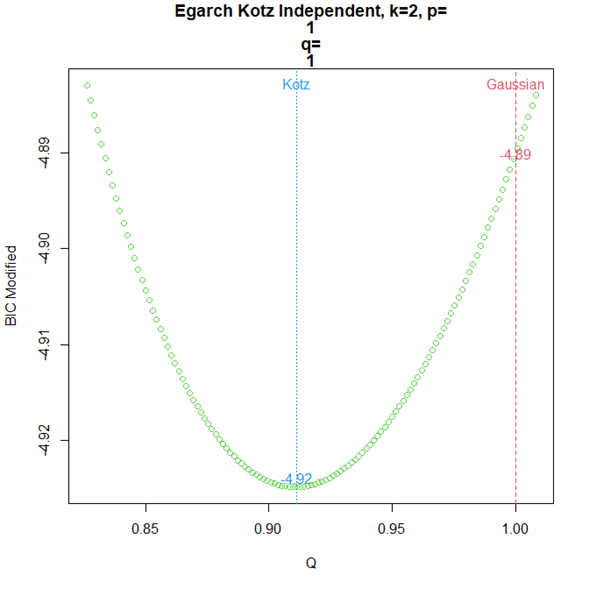}
         \caption{Yang and Yang profile criterion for $Q$.}
         \label{EgarchKotzk2p1q1f2}
     \end{subfigure}
        \caption{Egarch Kotz Independent ML estimation for $p=1, q=1$.}
       \label{EgarchKotzIndependent}
\end{figure}

Finally, the Egarch($p=1, q=1$) for an independent Kotz model provides the profile loglikelihood of parameter $Q$ in Figure \ref{EgarchKotzk2p1q1f1}, which can be set into a parallelism with the comparison by the modified BIC criterion of Yang and Yang of Figure \ref{EgarchKotzk2p1q1f2}.

\begin{figure}
         \centering
         \includegraphics[width=0.5\textwidth]{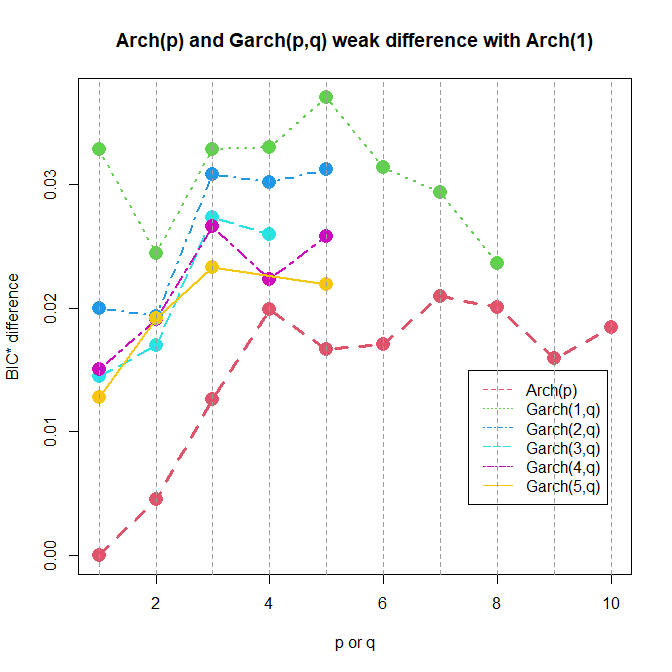}
         \caption{Arch(p) and Garch(p,q) BIC$^{*}$ difference with the Arch(1) .}
         \label{ArchGarchdifferenceWithArch1}
\end{figure}

Now, a common practice since the apparition of time series volatility laws ask for the best model of two o more models by means of an information criterion. AIC, BIC, HQC are the classical decision criteria and it can be checked in published real examples that the favoritism is based on minimum discrepancies. For example in the current time series data, commercial software and papers attain the best model for extremely low differences around $10^{-4}$, without a strong proof of the meaningful difference. The models here compared, usually are preferred in applications according the complexity is increased, then, in the Gaussian case, a Garch model tends to be better than Arch in some thousands of a decimal BIC, but, the Garch is worst in a similar small BIC amount than Tgarch, reports published no matter that the last models cannot be compared by BIC, because they are not nested. Then up here, we want to stress if each profile for $Q$ based on the new Yang and Yang criterion \citet{YangYang2007} of Figures \ref{ArchKotzk2p1f2}, \ref{GarchKotzk2p1q1f2}, \ref{TgarchKotzk2p1q1f2} and \ref{EgarchKotzk2p1q1f2} are showing that the respective Kotz model is better than the Gaussian case $Q=1$. It is not valid for this research, that the that apparently big differences in Yang and Yang criterion of 0.1, 0.06, 0.05 and 0.03, decides for a Kotz model of $Q=0.88, 0.91, 0.91, 0.91$ over a Gaussian model with $Q=1$, respectively, see Figures Figures \ref{ArchKotzk2p1f1}, \ref{GarchKotzk2p1q1f1}, \ref{TgarchKotzk2p1q1f1} and \ref{EgarchKotzk2p1q1f1}. 
In order to solve a robust selection we consider the comparison emerging from the Yang and Yang criterion \citet{YangYang2007} with the modified BIC$^{*}$ defined in (\ref{BIC*}). The corresponding comparison of the two models, follows \citet{KassRaftery1995} and \citet{Raftery1995} with the proposed grades of evidence in Table \ref{Table:BIC*} for BIC$^{*}$ difference.
Given that the Arch(p) and Garch(p,q) are nested models, we can compare them with the simple law. In this case we find the  BIC$^{*}$ difference with the Arch(1).  The results for some particular values of p and q can be found in Figure \ref{ArchGarchdifferenceWithArch1}.

A very weak difference evidence between all the models and Arch(1) is obtained. No matter how complex is the Garch(p,q), the time series data cannot reach a plausible value near for the minimum BIC$^{*}$ difference of 2. In contrast, the classical literature should refer Garch(5,5), in this example, as a robust model than the simplest one of Arch(1). The usual decision is taken from magnitud differences of $10^{-3}$ in AIC or BIC, HQC, but without passing the comparison by an evidence of the addressed research of \citet{KassRaftery1995} and \citet{Raftery1995}. At least in this data, which attains all the requirements for studying by Arch and Garch model, reveals no significant improvement of the complexity in the volatility description. For this counterexample, application of \citet{KassRaftery1995} and \citet{Raftery1995} in the time series setting opens and interesting question about the strong evidence of robustness in the variance models. 

\subsection{Gillette time series returns under dependent probabilistic samples based on the Pearson Type VII distribution}\label{subsec:GuilletteDependentPearsonVII}

We devote few lines to an application that considers the most important paradigm of this paper: the role of estimation on time series models with probabilistic independent samples, a fact that is considered in the core itself of the likelihood. Likelihood has been thought as an immutable concept in time series, but its philosophy fights strongly with the natural concept of a correlated random variable flexing in time. No much conceptual effort should be made for explaining that time series are extremely correlated, then its estimation, which is resumed in the likelihood, definitely should contemplate dependent probabilistic samples. In that sense, the likelihood reflecting the expected inner dynamics of the variable in time, naturally arrives to a realistic estimation of the volatility. This paradigm also appears in copula theory, which trying to solved the addressed evolution correlation, forgets the underlying and hidden basement of the likelihood and uses it for estimation of the parameters, buckling in the original intention of avoiding independence. Thus, implementation of a new perspective of the likelihood by a vector variate strongly dependent on time, should tackle the problem of volatility. We have seen the above mentioned effect in similar stages of financial models. It was the case of \citet{dgclpr:22}, where a classical likelihood decision based on the mean of a bivariate distribution was really placed in the tails of the right distribution estimated by a probabilistic dependent likelihood.  
Considering the again Gillette time series returns, if we provide a Pearson Type VII model for the underlying errors indexed by the non-estimated degrees of freedom $r$, then a plausible likelihood emerges in \ref{loglikelihoodPearsonDependent}. 
Inspired in the problem of very weak evidence for this special set of returns, and for the sake of illustration, we just estimate the simplest model of volatility Arch(1), Garch(1,1), Tgarch(1,1) and Egarch(1,1). As in the independent case, we use 8 methods of optimization, with gradient and Hessian support. All of them converges to the same estimation. Table \ref{Table:DependentModels} provides the results of each model and the corresponding comparison with the BIC$^{*}$ criterion for Arch and Garch and the following CAIC criterion for the remaining nested models (\citet{Bozdogan1987}).
\begin{equation}
    CAIC=\left[-2\mathcal{M}+n_{p}(\log(n)+1)\right]/n,
\end{equation}
where $\mathcal{M}$ is the maximum log-likelihood, $n$ is the sample size and $n_{p}$ is the number of estimated parameters. 
    \begin{table}[h!]
\begin{scriptsize}
\centering
\begin{tabular}{||c c c c c c c ||} 
 \hline
 Model&Shape parameter & BIC$^{*}$&CAIC&BIC$^{*}$ D.&CAIC D.&BIC$^{*}$ evidence\\ [0.5ex] 
 \hline\hline  
  Kotz Arch(1)& $Q=0.882733$& -4.86738&-4.844434&0&&\\
 Kotz Tgarch(1,1)& $Q=0.907778$&-4.925434&-4.894839&&0.050405&\\
 Kotz Egarch(1,1)&$Q=0.911211$ & -4.923831&-4.893235&&0.048801&\\
 \hline\hline
  Pearson Arch(2)& $r=1$&-4.829144&-4.806198&0&&\\
Pearson Arch(2)& $r=2$&-4.829654&-4.806708&0.00051&&Weak\\
Pearson Arch(2)& $r=3$&-4.829953& -4.807007&0.000809&&Weak\\
Pearson Arch(2)& $r=4$&-4.830171&-4.807225&0.001027&&Weak\\
Pearson Arch(2)& $r=5$&-4.830345&-4.807398&0.0012&&Weak\\
Pearson Garch(1,1)& $r=1$&-4.859952&-4.837006&0.030808&&Weak\\
Pearson Tgarch(1,1)& $r=1$&-4.895444&-4.868674&&0.062476&\\
Pearson Egarch(1,1)& $r=1$&-4.894745&-4.867974&&0.061776&\\ [1ex] 
 \hline
\end{tabular}
\caption{Comparison of independent and dependent models.}
\label{Table:DependentModels}
\end{scriptsize}
\end{table}

Table \ref{Table:DependentModels} refers the comparison of Kotz Tgarch(1,1) and Kotz Egarch(1,1) with the simplest model Kotz Arch(1) for independent probabilistic samples. Given that the three models are not nested, then we can only used the CAIC difference, but no evidence criteria for that magnitudes are available. This implies that no conclusion about the best model can be provided. The usual selection in favor of the best difference does not hold, i.e. we cannot affirm that Kotz Tgarch(1,1) is the best model. A similar not nested comparison is achieved in the dependent case for Pearson Egarch(1,1) and Pearson Tgarch(1,1), namely, we cannot assure that the  CAIC difference of 0.062476 defines the Pearson Tgarch(1,1) model over the remaining two laws.
However, when the models are nested, we can use \citet{KassRaftery1995} and \citet{Raftery1995} for an evidence of BIC$^{*}$ difference. This is the case of Table \ref{Table:DependentModels} when Pearson Garch(1,1) and Pearson Arch(2) with degrees of freedom $r=2,3,4,5$ are compared with the simplest model Pearson Arch(2) of 1 degree of freedom. A notorious weak evidence is obtained, and the BIC$^{*}$ differences are very far from the limit 2 for at least a positive predominance of Garch(1,1). Then the referred 5 models with increased complexity do not reach a best performance than the simplest Pearson Arch(2). In this sense, the rigorous evidence limits for a best model according to \citet{KassRaftery1995} and \citet{Raftery1995} opens in this example an interesting future work. Moreover, an further study of \citet{KassRaftery1995} and \citet{Raftery1995} into the setting of non nested models via CAIC should provide similar bounds of significant evidence. This important issue will be part of a future work.

\section{Conclusions}

This paper has proposed the so-called dependent time series, from a new perspective of likelihood based on dependent probabilistic samples governed by a general class of distributions that also allows multiple forms for volatility. The article also has taken from other areas the inclusion of robust selection criteria for classic volatility models and the requirement of degrees of evidence in the differentiation of the modified BIC criterion. The degrees of evidence that advance from weak to very strong, leave no opportunity for significant difference in a well-known time series, usually treated with hierarchical models that attribute a better explanation to the most complex model, but that under the incorporation of the new scale of evidence, does not achieve the minimum difference with the simplest model. This paradigm of dependent time series, generalized under families of distributions and hierarchized by demanding degrees of evidence, consolidates an alternative for artificial intelligence approaches and other increasingly complex volatility models, but with equal difficulty in activation functions or differentiation of the effect of volatility. Finally, given the small sample and non-seasonality events that strongly affect the foundations of classical volatility models, there are parallels to study in the future and to relate the new view of likelihood presented here with volatility modeled from the recent spectroscopic shape theory by \citet{VillarrealEtAl2019}.

\titleformat{\section}[block]{\normalfont\center #1}{}{0pt}{}

\section{Acknowledgements}
F.O.P. thanks the PhD in Modeling and Scientific Computing and the High Level Education Committee of the University of Medellín for supporting his doctoral studies. G.G. thanks the partial support of CBF2023-2024-3976, MX.

\end{document}